\documentclass[12pt,reqno]{amsart}

\newcommand{\SL}{\operatorname{SL}}
\newcommand\be{\begin{equation}}
\newcommand\ee{\end{equation}}

\usepackage[english,french]{babel}

\newtheorem*{thme}{Th\'eor\`eme}

\newtheorem*{thmnn}{Theorem}

\begin{document}

\author{Jean Bourgain}
\thanks{Bourgain is partially supported by NSF grant DMS-0808042.}
\email{bourgain@math.ias.edu}
\address{IAS, Princeton, NJ 08540}
\author{Alex Kontorovich}
\thanks{Kontorovich is partially supported by  NSF grants DMS-1064214 and DMS-1001252.}
\email{alexk@math.sunysb.edu}
\address{Stony Brook University, Stony Brook, NY 11794}

\title[ON ZAREBMA'S CONJECTURE]{ON ZAREBMA'S CONJECTURE\
\\
\
\\
SUR UNE CONJECTURE DE ZAREMBA
}

\selectlanguage{english}%

\date{\today}
\maketitle

\begin{abstract}
 It is shown that there is a constant $A$ and a density one subset
$S$ of the positive integers, such that for all $q\in S$ there is some $1\leq
p<q$, $(p, q)=1$, so that $\frac pq$ has all its partial quotients bounded by $A$.
\end{abstract}

 \selectlanguage{french}%
\begin{abstract}
On montre qu'il existe une constante $A$ et un
sous-ensemble $S$ des entiers positifs de densit\'e un, tel que pour tout
$q\in S$ il y a un entier $1\leq p<q, (p, q)=1$ pour lequel les quotients
partiels de $\frac pq$ sont born\'es par $A$.
\end{abstract}

\section*{Version Fran\c caise Abr\'egi\'ee}

Une conjecture due \`a Zaremba \cite{Zaremba1972} affirme qu'il existe pour tout entier $q\in
\Bbb Z_+$ un entier $1\leq p<q$, $(p, q)= 1$, tel que le d\'eveloppement $\frac pq =[a_1, \ldots, a_k]$ en fraction continue a tous ses quotients partiels $a_j$ born\'es
par une constante absolute $A$.
En fait, il est conjectur\'e que $A=5$ suffit. Notre r\'esultat principal est le suivant.

%
\begin{thme}
Il existe un sous-ensemble $S$ de $\Bbb Z_+$, de pleine densit\'e, v\'erifyant la conjecture de Zaremba pour une constante absolute $A$.
\end{thme}

Notre approche au probl\`eme est une adaptation de la m\'ethode introduite dans 
\cite{BourgainKontorovich2010} pour \'etudier les ensembles d'entiers g\'en\'er\'es par les orbites d'un sous-groupe
de $\SL_2(\Bbb Z)$.
La diff\'erence principale est que dans la situation pr\'esente, il s'agit d'un semi-groupe.
En particulier, nous ne faisons pas usage des r\'esultats de \cite{BourgainKontorovichSarnak2010}, mais plut\^ot de ceux de \cite{BourgainGamburdSarnak2009} bas\'es sur l'approche symbolique.

 \selectlanguage{english}%

{\bf 1.}
For given $A>0$, let $\mathcal C_A\subset [0, 1]$ be the Cantor-like set of real numbers $x$ in the unit interval, whose partial quotients are bounded by $A$.
Thus, writing $x$ in its contnued fraction expansion
$$
x=
\cfrac{1}{a_{1}+\cfrac{1}{a_{2}+\ddots\cfrac{1}{a_{k}+\ddots}}}
=
[a_1, a_2, \ldots, a_k, \ldots]
$$
we have that all partial quotients $a_k$ are bounded by $A$.
The Hausdorff dimension $\delta_A$ of $\mathcal C_A$ is asymptotically
$$
\delta_A =1-\frac 6{\pi^2A} - \frac {72 \log A}{\pi^4A^2}+ O\Big(\frac 1{A^2}\Big)
$$
as $A\to \infty$ \cite{Hensley1992}.

Let further $\mathcal R_A$ denote the set of all partial convergents $\frac pq, (p, q)=1$ of numbers in $\mathcal C_A$ and let $\mathcal Q_A$ be the set of all continuants $q$.

Zaremba's conjecture \cite{Zaremba1972} states that
$$
\mathcal Q_A=\Bbb Z_+
$$
for sufficiently large $A$ (possibly $A=5$).

It should be noted that the original motivation for this problem has to do with the theory of ``good'' lattice points and low-discrepancy  sequences in numerical
multi-dimensional integration and in pseudo-randomness.

It was shown by Niederreiter \cite{Niederreiter1986} that Zaremba's conjecture holds for small powers, in fact
$$
\{2^j\}, \{3^k\}
\subset\mathcal Q_3.
$$
On the other hand, a result due to Hensley \cite{Hensley1989} states that
$$
 N^{2\delta_A}\ll \# \Big\{\frac pq \in R_A; (p, q)=1  \ \text { and } \ 1\leq p< q\leq N\Big\} \ll N^{2\delta_A},
$$
an immediate consequence of which 
is that
$$
\# \mathcal Q_A\cap[1,N] \gg N^{2\delta_A-1}.
$$

Our main result is the following

\begin{thmnn} 
For $A$ sufficiently large, almost every integer satisfies Zaremba's conjecture. That is,
$$
\# \mathcal Q_A\cap[1,N]=N\big(1+o(1)\big)
,
$$
as $N\to\infty$.
($A=2189$ satisfies the claim).
\end{thmnn}

{\bf 2.}
A few comments about the method.

Our approach is an adaptation of the technique introduced in \cite{BourgainKontorovich2010} in the study of sequences of integers produced by orbits of subgroups $\Gamma$ of $\SL_2(\Bbb Z)$,
assuming the dimension $0<\delta<1$ of the limit set of $\Gamma$ close enough to 1.
We proceed by the Hardy-Littlewood circle method, analyzing certain relevant exponential sums on `minor' and `major' arcs.
While this approach is quite standard in number theoretical problems (for instance in the 
Goldbach problem), the ingredients involved in our situation are special.

In \cite{BourgainKontorovich2010}, the analysis on the minor arcs is achieved using Vinogradov-type multi-linear estimates, depending essentially on the group structure.
Then a precise evaluation of the exponential sum on the major arcs is obtained by 
relying on the spectral and representation theory of $\Gamma\backslash \SL_2$,
 as developed in \cite{BourgainKontorovichSarnak2010}.
The outcome is the usual local-to-global representation formula, with a small exceptional set.

It turns out that Zaremba's conjecture admits a formulation of similar flavor.
Let $\mathcal G_A$ be the semi-group generated by the matrices
$$
\begin{pmatrix} 0&1\\ 1&a\end{pmatrix}  \ \text { with } \ 1\leq a\leq A
$$
and observe that
$$
\begin{pmatrix} 0&1\\1&a_1\end{pmatrix} \begin{pmatrix} 0&1\\ 1&a_2\end{pmatrix}
\cdots \begin{pmatrix} 0&1\\ 1&a_k\end{pmatrix} =\begin{pmatrix} *&p\\ *&q\end{pmatrix}
$$
is 
equivalent to $\frac pq =[a_1, \ldots, a_k]$.

Hence, the orbit $\mathcal G_Ae_2$, with $e_2 =(0, 1)$, consists precisely of the set of coprime pairs $(p, q)$ with $\frac pq\in R_A$.

The main difference with \cite{BourgainKontorovich2010} is that instead of the group $\Gamma$, the semi-group $\mathcal G_A$ is involved.
It turns out however that this distinction has essentially no effect on the minor arcs analysis.
On the other hand to proceed with the description of the exponential sum on the major arcs, the automorphic approach from \cite{BourgainKontorovichSarnak2010} is no longer applicable.
Instead we rely on the 
thermodynamical formalism
based on the Ruelle transfer operator (which actually is already exploited in
\cite{Hensley1989}).
Here our aim is to establish certain equidistributional properties from a joint  Archimedean / modular perspective, since this  allows us to analyze the exponential sum on a
major arc $\theta =\frac rs+\beta$, $(r, s)=1$ with $s<N^\epsilon$ and $|\beta|<1/N^{1-\epsilon}$.
This type of (quantitative) result is provided by \cite{BourgainGamburdSarnak2009} in a form applicable to our problem.

\bibliographystyle{alpha}
\bibliography{../../AKbibliog}
\end{document}